\documentclass[11pt,reqno]{amsart}
\usepackage{graphicx}
\usepackage{amsmath, amssymb}
\usepackage{verbatim}
\usepackage{color}

\newtheorem{conj.}[thm]{Conjecture}

\theoremstyle{definition}

\theoremstyle{remark}

\numberwithin{equation}{section}

\begin{document}

\title[Uncertainty Principles for Quaternion Windowed Offset Linear Canonical Transform of Two Dimensional Signals]
{Uncertainty Principles for Quaternion Windowed Offset Linear Canonical Transform of Two Dimensional Signals}

\author[Aajaz A. Teali]{Aajaz A. Teali}

\address{ Department of Mathematics, University of Kashmir, South Campus, Anantnag-192101 Jammu and Kashmir, India.}
\email{aajaz.math@gmail.com}


\begin{abstract}
The offset linear canonical transform encompassing the numerous integral transforms, is a promising tool for analyzing non-stationary signals with more degrees of freedom. In this paper, we generalize the windowed offset linear canonical transform for quaternion-valued signals, by introducing a novel time-frequency transform namely the quaternion windowed offset linear canonical transform of 2D quaternion-valued signals. We initiate our investigation by studying some fundamental properties of the proposed transform including inner product relation, energy conservation, and reproducing formula by employing the machinery of quaternion offset linear canonical transforms. Some uncertainty principles such as Heisenberg-Weyl, logarithmic and local uncertainty principle are also derived for quaternion windowed offset linear canonical transform. Finally, we gave an example of quaternion windowed offset linear canonical transform.
\end{abstract}

\parindent=0mm \vspace{.0in}
\maketitle
\section{Introduction}

\parindent=0mm \vspace{.0in}
During the culminating years of last century, a reasonably flexible integral transform associated with a general inhomogeneous lossless linear mapping in phase-space was introduced, namely, the offset linear canonical transform(OLCT)  \cite{AS,Cai}. The OLCT is a six-parameter class of linear integral transform which encompasses a number of well known unitary transforms including the classical Fourier transform, fractional Fourier transform, Fresnel transform, Laplace transform, Gauss-Weierstrass transform, Bargmann transform and the linear canonical transform \cite{DS1,Al,HKOS}. Due to the extra degrees of freedom, OLCT has attained a respectable status within a short span and is being broadly employed across several disciplines of science and engineering including signal and image processing, optical and radar systems, electrical and communication systems, pattern recognition, sampling theory, shift-invariant theory and quantum mechanics \cite{Qzy,Z,Xhch,ULT,BZ}. Recently, in \cite{Shah} we introduce a hybrid integral transform namely, windowed special affine Fourier transform which is capable of providing a joint time and frequency localization of non-stationary signals with more degrees of freedom. In view of its numerous applications, one is particularly interested in its higher dimensional analogues.

\parindent=8mm \vspace{.1in}
In the meantime, quaternion algebra has become an active area of research as it offers a refinement of classical harmonic analysis and is used to generalize the classical theory of holomorphic functions of one complex variable onto the multidimensional situation. It gives a simple and profound representation of signals wherein several components are to be controlled simultaneously. The development of integral transforms for quaternion valued signals has found numerous applications in 3D computer graphics, aerospace engineering, artificial intelligence and colour image processing. The extension of classical Fourier transform to quaternion algebra has been introduced in \cite{Ell} and its efficient implication was given by Pei et al. in \cite{Ding}. The QFTs play a vital role in the representation of quaternion signals. They transform a 2D quaternion-valued signals into the quaternion-valued frequency domain signals. Many applications of the QFTs have been found in color image processing, especially in color-sensitive smoothing, speech recognition, edge detection and data compression \cite{Bas,Bayro,Dubey,Gri}. Moreover, some authors have also generalized the classical linear canonical transform to quaternion-valued signals, known as the quaternion linear canonical transform (QLCT). The QLCT was firstly studied in \cite{Kou2} including prolate spheroidal wave signals and uncertainty principles \cite{Kou}. Some useful properties and applications of the QLCT such as linearity, reconstruction formula, continuity, boundedness, positivity inversion formula and the uncertainty principle were established in \cite{Asym,MR,Zhang} and \cite{BahApp}. Recently, based on the (two-sided) QLCT \cite{Kou2}, the quaternion windowed linear canonical transform of 2D quaternion signals has been introduced by Gao and Li in \cite{Wen}, which generalizes the quaternion windowed Fourier transform \cite{MawW}.  In \cite{QOLCT}, Haoui et al. introduced the quaternion offset linear canonical transforms(QOLCT) as a generalization of QLCTs and QFTs. But it fails in obtaining the local features of non-transient signals due to its global kernel, it is therefore interesting and worthwhile to investigate the localization of quaternion offset linear canonical transforms.

\parindent=8mm \vspace{.1in}
In the present work, we study the generalization of the windowed offset linear canonical transform for quaternion-valued signals, we introduce a novel time-frequency transform namely the quaternion windowed offset linear canonical transform of 2D quaternion-valued signals. It can reveal the local QOLCT-frequency contents and enjoys high concentrations and eliminates the cross term. Some important properties are analyzed. Moreover, inner product relation, energy conservation and inversion formula are establised. Furthermore, some uncertainty principles such as Heisenberg-Weyl, logarithmic and local uncertainty principle are also derived for quaternion windowed offset linear canonical transform. Finally, we gave an example of quaternion windowed offset linear canonical transform.

\parindent=8mm \vspace{.1in}
The organization of the article is as follows: We begin in Section 2 by presenting the brief recall of quaternion algebra and different types of quaternion offset linear canonical transforms. In Section 3, we introduce the concept of quaternion windowed offset linear canonical transform and obtain the expected properties of the QWOLCT including Parseval's formula, energy conservation, isometry and inversion formula. Moreover, the uncertainty principles such as Heisenberg-Pauli-Weyl inequality, logarithmic uncertainty principle and local uncertainty principle are generalized in the quaternion windowed offset linear canonical domains in Section 4. Finally, an example of quaternion windowed offset linear canonical transform is given in Section 5.

\section{Quaternion Algebra and Quaternion Offset Linear Canonical Transform}

\subsection{ Quaternion Algebra}
\parindent=0mm \vspace{.1in}
The theory of quaternions was initiated by the Irish mathematician Sir W.R. Hamilton in 1843 and is denoted by $\mathbb H$ in his honour.
The quaternion algebra provides an extension of the complex number system to an associative non-commutative four-dimensional algebra.
The quaternion algebra $\mathbb H$ over $\mathbb R$ is given by
\begin{align*}
\mathbb H=\Big\{f=a_{0}+i\,a_{1}+j\,a_{2}+k\,a_{3}\,:\,a_{0},a_{1},a_{2},a_{3}\in\mathbb R\Big\},
\end{align*}
where $i,j,k$ denote the three imaginary units, obeying the Hamilton's multiplication rules
\begin{align*}
ij=k=-ji,~jk=i=-kj\,,\,ki=j=-ik\,,\, {\text {and}}\,\, {i}^{2}={j}^{2}={k}^{2}={ijk}=-1.
\end{align*}
For  quaternions $f_{1}=a_{0}+i\,a_{1}+j\,a_{2}+k\,a_{3}$ and $f_{2}=b_{0}+i\,b_{1}+ j\,b_{2}+k\,b_{3}$, the addition  is defined componentwise
and the multiplication is defined as
\begin{align*}
f_{1}f_{2}&=(a_{0}b_{0}-a_{1}b_{1}-a_{2}b_{2}-a_{3}b_{3})+i\,(a_{1}b_{0}+a_{0}b_{1}+a_{2}b_{3}-a_{3}b_{2})\\
&~~~~~~~~~~~~~~+ j(a_{0}b_{2}+a_{2}b_{0}+a_{3}b_{1}-a_{1}b_{3})+k\,(a_{0}b_{3}+a_{3}b_{0}+a_{1}b_{2}-a_{2}b_{1}).
\end{align*}
The conjugate and norm of a quaternion  $f=a_{0}+i\,a_{1}+j\,a_{2}+k\,a_{3},$ are given by $\overline{f}=a_{0}-i\,a_{1}-j\,a_{2}-k\,a_{3}$ and ${\|f\|}_{\mathbb H}=\sqrt{{a_{0}}^{2}+{a_{1}}^{2}+{a_{2}}^{2}+{a_{3}}^{2}}$, respectively. We also note that an  arbitrary quaternion $f$ can be represented by two complex numbers  as $f=(a_{0}+i\,a_{1})+j\,(a_{2}-i\,a_{3})=f_1+j\,f_2$, where $f_1,f_2\in\mathbb C$, and hence, ${\overline{f}}={\overline f_1}-j\,f_2$, with $\overline{f_1}$ denoting the complex conjugate of $f_1$. Moreover, the inner product of any two quaternions  $f=f_{1}+\,jf_{2},$ and $g=g_{1}+j\,g_{2}$ in ${\mathbb H}$ is defined by
\begin{align*}
{\big\langle f,g\big\rangle}_{\mathbb H}&= f\overline{g}=(f_{1}{\overline g_{1}}+{\overline f_{2}}g_{2})+j(f_{2}{\overline g_{1}}-{\overline f_{1}}g_{2}).
\end{align*}

By virtue of the complex domain representation, a quaternion-valued function $f:\mathbb R^2 \to \mathbb H$ can be decomposed as $f(x)=f_{1}+j\,f_{2}$,  where $f_{1},f_{2}$ are both complex valued functions.

\parindent=8mm \vspace{.1in}
Let us denote $L^2(\mathbb R^2,\mathbb H)$, the space of all quaternion valued functions $f$ satisfying
\begin{align*}
 {\big\|f\big\|}_{2}=\left\{\int_{\mathbb R^{2}}\Big(|f_{1}(x)|^{2}+|f_{2}(x)|^{2}\Big)\,dx\right\}^{1/2}< \infty.
\end{align*}
The norm on $L^2(\mathbb R^2,\mathbb H)$ is obtained from the inner product of the quaternion valued functions $f=f_{1}+j\,f_{2},$ and $g=g_{1}+j\,g_{2}$  as
\begin{align*}
{\big\langle f,\,g\big\rangle}_{2}&= \int_{\mathbb R^{2}}{\big\langle f,\,g\big\rangle}_{\mathbb H}\,dx\\
&=\int_{\mathbb R^{2}}\bigg\{\Big(f_{1}(x)\,\overline{g_{1}}(x)+\overline{f_{2}}(x)\,g_{2}(x)\Big)+j \Big(f_{2}(x)\,\overline{g_{1}}(x)
-\overline{f_{1}}(x)\,g_{2}(x)\Big)\bigg\}\,dx.
\end{align*}
An easy computation shows that $L^2(\mathbb R^2,\mathbb H)$ equipped with above defined inner product is a Hilbert space.

\parindent=8mm \vspace{.1in}
Before we proceed to establish the main results of this paper, we shall recall the basic theory of quaternion offset linear canonical transform

\subsection{ Quaternion Offset Linear Canonical Transform }
Due to non-commutativity of quaternion multiplication, there are three types of the quaternion offset linear canonical transform, the left-sided QOLCT, the right-sided QOLCT, and two-sided QOLCT.  In this paper, we will extend the theory of two-sided QOLCT introduced by Haoui et al. \cite{QOLCT} as follows:

\parindent=0mm \vspace{.1in}
{\bf Definition 2.1. } Let $A_s =\left[\begin{array}{cccc}a_s & b_s &| & p_s\\c_s & d_s &| & q_s \\\end{array}\right] $, be a matrix parameter such that $a_s$, $b_s$, $c_s$, $d_s$, $p_s$, $q_s \in \mathbb R$ and $ a_sd_s-b_sc_s=1,$ for $s=1,2.$ The two-sided quaternion offset linear canonical transform of any quaternion valued function $f\in L^2(\mathbb R^2,\mathbb H)$, is given by
\begin{align*}
\mathcal O_{A_1,A_2}^{i,j}\big[f\big](w)=\left\{\begin{array}{cc} \int_{\mathbb R^2} K_{A_1}^i(t_1,w_1) \,f(t)\,K_{A_2}^j(t_2,w_2)dt, &b_1,b_2 \neq 0,\\ \sqrt{d_1}\,e^{i\big[\frac{c_1d_1}{2}(w_1-p_1)^2+w_1p_1\big]} f(d_1w_1-d_1p_1,\,t_2) K_{A_2}^j(t_2,w_2), &b_1=0,b_2\neq0,\\
\sqrt{d_2}K_{A_1}^i(t_1,w_1) f(t_1,d_2w_2-d_2p_2) \,e^{j\big[\frac{c_2d_2}{2}(w_2-p_2)^2+w_2p_2\big]}, &b_1\neq0,b_2=0,\\
\sqrt{d_1d_2}f(d_1w_1-d_1p_1,d_2w_2-d_2p_2)\,e^{i\big[\frac{c_1d_1}{2}(w_1-p_1)^2+w_1p_1\big]}& \\ \qquad\times\,e^{j\big[\frac{c_2d_2}{2}(w_2-p_2)^2+w_2p_2\big]}, &b_1=0,b_2=0,
\end{array}\right.
\end{align*}
where $t=(t_{1},t_{2}),\, w=(w_{1},w_{2})$ and the quaternion kernels $K_{A_1}^i(t_1,w_1)$ and $K_{A_2}^j(t_2,w_2)$ are respectively given by
\begin{align*}
K_{A_1}^i(t_1,w_1)=\dfrac{1}{\sqrt{2\pi b_1}} \,e^{\frac{i}{2b_1}\big[a_1t_1^2-2t_1(w_1-p_1)-2w_1(d_1p_1-b_1q_1)+d_1(w_1^2+p_1^2)-\frac{\pi b_1}{2}\big]},b_1\neq0\,\tag{2.1}
\end{align*}

\begin{align*}
K_{A_2}^j(t_2,w_2)= \dfrac{1}{\sqrt{2\pi b_2}} \,e^{\frac{j}{2b_2}\big[a_2t_2^2-2t_2(w_2-p_2)-2w_2(d_2p_2-b_2q_2)+d_2(w_2^2+p_2^2)-\frac{\pi b_2}{2}\big]},\,b_2\neq0\,\tag{2.2}
\end{align*}

{\bf Remark 2.2.} The left-sided and right-sided QOLCT can be defined correspondingly by placing the two above kernels both on the left or on the right, respectively.

The corresponding inversion formula for two-sided QOLCT is given by
\begin{align*}
f(x) &=\int_{\mathbb R^2}\overline{K_{A_1}^i(t_1,w_1)} \, \mathcal O_{A_1,A_2}^{i,j}\big[f\big](w)\,\overline{K_{A_2}^j(t_2,w_2)} \,d\xi.\tag{2.3}
\end{align*}

\parindent=0mm \vspace{.1in}
{\bf Theorem 2.3 (Plancherel for QOLCT).} {\it Every two dimensional quaternion valued function $f \in L^2(\mathbb R^2,\mathbb H)$ and its two-sided QOLCT are related to the Plancherel identity in the following way:}
\begin{align*}
\big\|\mathcal O_{A_1,A_2}^{i,j}\big[f\big]\big\|_{L^2(\mathbb R^2,\mathbb H)}=\left\|f\right\|_{L^2(\mathbb R^2,\mathbb H)}.\,\tag{2.4}
\end{align*}

\parindent=0mm \vspace{.0in}

\section{  Quaternion Windowed Offset Linear Canonical Transform}

In this section, we shall formally introduce the notion of the two-sided quaternion windowed offset linear canonical transform of two dimensional signals and then establish the fundamental properties of the proposed transform.

\parindent=0mm \vspace{.1in}
{\bf Definition 3.1.} Let $A_s =\left[\begin{array}{cccc}a_s & b_s &| & p_s\\c_s & d_s &| & q_s \\\end{array}\right] $, be a matrix parameter such that $a_s$, $b_s$, $c_s$, $d_s$, $p_s$, $q_s \in \mathbb R$ and $ a_sd_s-b_sc_s=1,$ for $s=1,2.$ The two-sided quaternion windowed offset linear canonical transform of any quaternion valued function $f\in L^2(\mathbb R^2,\mathbb H)$, with respect to a non-zero window function $g\in  L^2(\mathbb R^2,\mathbb H)$ is given by  \begin{align*}\label{3.1}
\mathcal O_{g}^{A_1,A_2}\big[f\big](u,w)=\left\{\begin{array}{cc} \int_{\mathbb R^2} K_{A_1}^i(t_1,w_1) \,f(t)\overline{g(t-u)}\,K_{A_2}^j(t_2,w_2)dt, &b_1,b_2 \neq 0,\\ \sqrt{d_1} \,e^{i\big[\frac{c_1d_1}{2}(w_1-p_1)^2+w_1p_1\big]} f(d_1w_1-d_1p_1,t_2) K_{A_2}^j(t_2,w_2), &b_1=0,b_2\neq0,\\
\sqrt{d_2} K_{A_1}^i(t_1,w_1) f(t_1,d_2w_2-d_2p_2) \quad\qquad\qquad\qquad\qquad & \\ \times\overline{g(t_1-u_1,d_2w_2-d_2p_2)}\,e^{j\big[\frac{c_2d_2}{2}(w_2-p_2)^2+w_2p_2\big]}, &b_1\neq0,b_2=0,\\
\sqrt{d_1d_2} f(d_1w_1-d_1p_1,d_2w_2-d_2p_2)\quad\qquad\qquad\qquad\qquad& \\ \times\overline{g (d_1w_1-d_1p_1,d_2w_2-d_2p_2) }\,e^{i\big[\frac{c_1d_1}{2}(w_1-p_1)^2+w_1p_1\big]}& \\ \times e^{j\big[\frac{c_2d_2}{2}(w_2-p_2)^2+w_2p_2\big]}, &b_1=0,b_2=0,
\end{array}\right.\,\tag{3.1}
\end{align*}
where $t=(t_{1},t_{2}),\, w=(w_{1},w_{2}),\, u= (u_1,u_2)$ and the quaternion kernels $K_{A_1}^i(t_1,w_1)$ and $K_{A_2}^j(t_2,w_2)$ are given by (2.1) and (2.2), respectively.

\parindent=0mm \vspace{.1in}
{\bf Remark:} It is worth to note that the quaternion windowed OLCT (3.1), boils down to various linear integral transforms such as:
\begin{itemize}
\item Windowed versions of quaternion linear canonical transform when matrices parameters $A_s =\left[\begin{array}{cccc}a_s & b_s &| & 0\\c_s & d_s &| & 0 \\\end{array}\right] $,\\

\item Quaternion windowed fractional Fourier transform when $A_s =\left[\begin{array}{cccc}\cos\theta & \sin\theta &| & 0\\-\sin\theta & \cos\theta &| & 0 \\\end{array}\right] $,\\

\item Quaternion windowed Fourier transform when $A_s =\left[\begin{array}{cccc}1 & 0 &| & 0\\0 & 1 &| & 0 \\\end{array}\right]$.

\end{itemize}
\parindent=0mm \vspace{.1in}
For brevity, in this paper we focus for the case $b_s\neq 0,\, s=1,2,$ as in other cases proposed transform reduces to a chrip multiplications.

\parindent=0mm \vspace{.1in}
In the sequel, our intention is to study the fundamental properties of the proposed quaternion windowed OLCT (3.1).

\parindent=0mm \vspace{.1in}
{\bf Property 1.} {\it (Linearity) For any given quaternion function $f_n \in L^2(\mathbb R^2,\mathbb H), n\in \mathbb N,$ the following relationship is true:}
\begin{align*}
\mathcal O_{g}^{A_1,A_2}\Big[\sum_{n\in \mathbb N}\alpha_n f_n\Big](u,w)= \sum_{n\in \mathbb N}\alpha_n \mathcal O_{g}^{A_1,A_2}\big[ f_n\big](u,w),\,\alpha_n\in \mathbb H.
\end{align*}
\parindent=0mm \vspace{.1in}
{\bf Proof.} Since quaternion windowed OLCT is one of the linear integral transforms, Property 1 is directly obtained form the definition 3.1.

\parindent=0mm \vspace{.1in}
{\bf Property 2.} {\it (Time-Shift) For any given quaternion function $f \in L^2(\mathbb R^2,\mathbb H),$ and window function $g \in L^2(\mathbb R^2,\mathbb H),$ we have:}
\begin{align*}
{\mathcal{O}}^{A_1,A_2}_g\Big[f(t-k)\Big](u,w)&=\,e^{i\big[a_1k_1q_1-c_1k_1p_1+c_1k_1w_1-\frac{a_1c_1k_1^2}{2}\big]}\, {\mathcal{O}}^{A_1,A_2}_g  \big[f(t)\big](u-k,w-ak)\\
&\qquad\qquad\qquad\times\,e^{j\big[a_2k_2q_2-c_2k_2p_2+c_2k_2w_2-\frac{a_2c_2k_2^2}{2}\big]},
\end{align*}
{\it where $u=(u_1,u_2),\,k=(k_1,k_2),\,w=(w_1,w_2),\, a=(a_1,a_2).$}

\parindent=0mm \vspace{.1in}
{\bf Proof.} By Definition 3.1, we have
\begin{align*}
&{\mathcal{O}}^{A_1,A_2}_g\big[f(t-k)\big](u,w)\\
&\qquad=\int_{\mathbb R^2} K_{A_1}^i(t_1,w_1) \,f(t-k)\overline{g(t-u)}\,K_{A_2}^j(t_2,w_2)dt\\
&\qquad=\int_{\mathbb R^2} K_{A_1}^i(x_1+k_1,w_1) \,f(x)\overline{g(x+k-u)}\,K_{A_2}^j(x_2+k,w_2)dx\\
&\qquad=\dfrac{1}{2\pi \sqrt{b_1b_2}}\int_{\mathbb R^2} \,e^{\frac{i}{2b_1}\big[a_1(x_1+k_1)^2-2(x_1+k_1)(w_1-p_1)-2w_1(d_1p_1-b_1q_1)+d_1(w_1^2+p_1^2)-\frac{\pi b_1}{2}\big]}\\
&\qquad\qquad\times f(x)\overline{g\big(x-(u-k)\big)}\\
&\qquad\qquad\times\,e^{\frac{j}{2b_2}\big[a_2(x_2+k_2)^2-2(x_2+k_2)(w_2-p_2)-2w_2(d_2p_2-b_2q_2)+d_2(w_2^2+p_2^2)-\frac{\pi b_2}{2}\big]}dx\\
&\qquad=\dfrac{1}{2\pi \sqrt{b_1b_2}}\int_{\mathbb R^2} \,e^{\frac{i}{2b_1}\big[a_1k_1^2-2k_1(w_1-p_1)-2a_1k_1(d_1p_1-b_1q_1)-d_1(a_1^2k_1^2+2w_1a_1k_1)\big]}\\
&\qquad\qquad\times\,e^{\frac{i}{2b_1}\big[a_1x_1^2-2x_1(w_1-p_1-a_1k_1)-2(w_1-a_1k_1)(d_1p_1-b_1q_1)+d_1\big((w_1-a_1k_1)^2+p_1^2\big)-\frac{\pi b_1}{2}\big]}\\
&\qquad\qquad\times f(x)\overline{g\big(x-(u-k)\big)}\\
&\qquad\qquad\times\,e^{\frac{j}{2b_2}\big[a_2x_2^2-2x_2(w_2-p_2-a_2k_2)-2(w_2-a_2k_2)(d_2p_2-b_2q_2)+d_2\big((w_2-a_2k_2)^2+p_2^2\big)-\frac{\pi b_2}{2}\big]}dx\\
&\qquad\qquad\times\,e^{\frac{j}{2b_2}\big[a_2k_2^2-2k_2(w_2-p_2)-2a_2k_2(d_2p_2-b_2q_2)-d_2(a_2^2k_2^2+2w_2a_2k_2)\big]}\\
&\qquad= \,e^{\frac{i}{2b_1}\big[a_1k_1^2-2k_1(w_1-p_1)-2a_1k_1(d_1p_1-b_1q_1)-d_1(a_1^2k_1^2+2w_1a_1k_1)\big]}\\
&\qquad\qquad\times\int_{\mathbb R^2} K_{A_1}^i(x_1,w_1-a_1k_1)f(x)\overline{g\big(x-(u-k)\big)}K_{A_2}^j(x_2,w_2-a_2k_2)\, dx\\
&\qquad\qquad\times\,e^{\frac{j}{2b_2}\big[a_2k_2^2-2k_2(w_2-p_2)-2a_2k_2(d_2p_2-b_2q_2)-d_2(a_2^2k_2^2+2w_2a_2k_2)\big]}\\
&\qquad=\,e^{i\big[a_1k_1q_1-c_1k_1p_1+c_1k_1w_1-\frac{a_1c_1k_1^2}{2}\big]}\, {\mathcal{O}}^{A_1,A_2}_g  \big[f(t)\big](u-k,w-ak)\\
&\qquad\qquad\qquad\qquad\qquad\qquad\qquad\qquad\times\,e^{j\big[a_2k_2q_2-c_2k_2p_2+c_2k_2w_2-\frac{a_2c_2k_2^2}{2}\big]},
\end{align*}
where we used the relation $ a_sd_s=1+b_sc_s,$ for $s=1,2,$ in getting last equality.

\parindent=0mm \vspace{.1in}
{\bf Property 3.} {\it (Modulation) For any given quaternion function $f \in L^2(\mathbb R^2,\mathbb H),$ and window function $g \in L^2(\mathbb R^2,\mathbb H),$ then if $\mathcal M_w f(t)= e^{it_1w_1}f(t)e^{jt_2w_2},$  we have:}
\begin{align*}
&{\mathcal{O}}^{A_1,A_2}_g\Big[\mathcal M_w f(t)\Big](u,w)\\
&\quad=\,e^{iw_1(b_1q_1-d_1p_1)-i\frac{d_1b_1w_1^2}{2}+id_1w_1^2}\, {\mathcal{O}}^{A_1,A_2}_g  \big[f(t)\big]\big(u,(1-b)w\big)\,e^{jw_2(b_2q_2-d_2p_2)-j\frac{d_2b_2w_2^2}{2}+jd_2w_2^2},
\end{align*}
{\it where $u=(u_1,u_2),\,w=(w_1,w_2),\, a=(a_1,a_2).$}

\parindent=0mm \vspace{.1in}
{\bf Proof.} We have

\begin{align*}
&{\mathcal{O}}^{A_1,A_2}_g\Big[\mathcal M_w f(t)\Big](u,w)\\
&\qquad= \int_{\mathbb R^2} K_{A_1}^i(t_1,w_1) \,e^{it_1w_1}f(t)e^{jt_2w_2}\,\overline{g(t-u)}\,K_{A_2}^j(t_2,w_2)dt\\
&\qquad=\dfrac{1}{2\pi \sqrt{b_1b_2}}\int_{\mathbb R^2} \,e^{\frac{i}{2b_1}\big[a_1x_1^2-2x_1(w_1-p_1)-2w_1(d_1p_1-b_1q_1)+d_1(w_1^2+p_1^2)-\frac{\pi b_1}{2}\big]}\\
&\qquad\qquad\qquad\times \,e^{it_1w_1}f(t)e^{jt_2w_2}\,\overline{g(t-u)}\,e^{\frac{j}{2b_2}\big[a_2x_2^2-2x_2(w_2-p_2)-2w_2(d_2p_2-b_2q_2)+d_2(w_2^2+p_2^2)-\frac{\pi b_2}{2}\big]}dx\\
&\qquad=\dfrac{1}{2\pi \sqrt{b_1b_2}}\int_{\mathbb R^2} \,e^{\frac{i}{2b_1}\big[a_1x_1^2-2x_1(w_1-b_1w_1-p_1)-2w_1(d_1p_1-b_1q_1)+d_1(w_1^2+p_1^2)-\frac{\pi b_1}{2}\big]}\\
&\qquad\qquad\qquad\qquad\times \,f(t)\overline{e^{-jt_2w_2}}\,\overline{g(t-u)}\,\overline{\,e^{\frac{-j}{2b_2}\big[a_2x_2^2-2x_2(w_2-p_2)-2w_2(d_2p_2-b_2q_2)+d_2(w_2^2+p_2^2)-\frac{\pi b_2}{2}\big]}}dx\\
&\qquad=\dfrac{1}{2\pi \sqrt{b_1b_2}}\int_{\mathbb R^2} \,e^{\frac{i}{2b_1}\big[a_1x_1^2-2x_1(w_1-b_1w_1-p_1)-2w_1(d_1p_1-b_1q_1)+d_1(w_1^2+p_1^2)-\frac{\pi b_1}{2}\big]}\\
&\qquad\qquad\qquad\qquad\times \,f(t)\overline{\,e^{\frac{-j}{2b_2}\big[a_2x_2^2-2x_2(w_2-b_2w_2-p_2)-2w_2(d_2p_2-b_2q_2)+d_2(w_2^2+p_2^2)-\frac{\pi b_2}{2}\big]}\,g(t-u)}dx\\
&\qquad=\dfrac{1}{2\pi \sqrt{b_1b_2}}\int_{\mathbb R^2} \,e^{\frac{i}{2b_1}\big[a_1x_1^2-2x_1(w_1-b_1w_1-p_1)-2w_1(d_1p_1-b_1q_1)+d_1(w_1^2+p_1^2)-\frac{\pi b_1}{2}\big]}\\
&\qquad\qquad\qquad\qquad\times \,f(t)\overline{g(t-u)}\,e^{\frac{j}{2b_2}\big[a_2x_2^2-2x_2(w_2-b_2w_2-p_2)-2w_2(d_2p_2-b_2q_2)+d_2(w_2^2+p_2^2)-\frac{\pi b_2}{2}\big]}\,dx\\
&\qquad=\dfrac{1}{2\pi \sqrt{b_1b_2}}\,e^{\frac{i}{2b_1}\big[-2b_1w_1(d_1p_1-b_1q_1)-d_1b_1^2w_1^2+2d_1b_1w_1^2\big]}\\
&\qquad\quad\times\int_{\mathbb R^2} \,e^{\frac{i}{2b_1}\big[a_1x_1^2-2x_1(w_1-b_1w_1-p_1)-2(w_1-b_1w_1)(d_1p_1-b_1q_1)+d_1\big((w_1-b_1w_1)^2+p_1^2\big)-\frac{\pi b_1}{2}\big]}\\
&\qquad\quad\times f(t)\overline{g(t-u)}\,e^{\frac{j}{2b_2}\big[a_2x_2^2-2x_2(w_2-b_2w_2-p_2)-2(w_2-b_2w_2)(d_2p_2-b_2q_2)+d_2\big((w_2-b_2w_2)^2+p_2^2\big)-\frac{\pi b_2}{2}\big]}\,dx\\
&\qquad\quad\times\,e^{\frac{j}{2b_2}\big[-2b_2w_2(d_2p_2-b_2q_2)-d_2b_2^2w_2^2+2d_2b_2w_2^2\big]}\\
&\qquad=\,e^{iw_1(b_1q_1-d_1p_1)-i\frac{d_1b_1w_1^2}{2}+id_1w_1^2}\int_{\mathbb R^2} K_{A_1}^i(t_1,w_1-b_1w_1)\,f(t)\overline{g(t-u)}\, K_{A_1}^j(t_2,w_2-b_2w_2)\,dx\\
&\qquad\qquad\qquad\qquad\qquad\qquad\qquad\qquad\qquad\qquad\times\,e^{jw_2(b_2q_2-d_2p_2)-j\frac{d_2b_2w_2^2}{2}+jd_2w_2^2}\\
&=\,e^{iw_1(b_1q_1-d_1p_1)-i\frac{d_1b_1w_1^2}{2}+id_1w_1^2}\, {\mathcal{O}}^{A_1,A_2}_g  \big[f(t)\big]\big(u,(1-b)w\big) \,e^{jw_2(b_2q_2-d_2p_2)-j\frac{d_2b_2w_2^2}{2}+jd_2w_2^2}.
\end{align*}

\parindent=0mm \vspace{.1in}
{\bf Property 4.} {\it (Parity) Let $g \in L^2(\mathbb R^2,\mathbb H),$ be a quaternion window function and $f \in L^2(\mathbb R^2,\mathbb H),$ then we have:}
\begin{align*}
{\mathcal{O}}^{A_1,A_2}_g\Big[P f(t)\Big](u,w)=-{\mathcal{O}}^{A_1^\prime,A_2^\prime}_{Pg}\big[f(t)\big](-u,-w),
\end{align*}
{\it where $Pf(t)=f(-t),$ and $ A_s^\prime =\left[\begin{array}{cccc}a_s & b_s &| & -p_s\\c_s & d_s &| & -q_s \\\end{array}\right],\,s=1,2. $}

\parindent=0mm \vspace{.1in}
{\bf Proof.} We have

\begin{align*}
&{\mathcal{O}}^{A_1,A_2}_g\Big[P f(t)\Big](u,w)\\
&\qquad= \int_{\mathbb R^2} K_{A_1}^i(t_1,w_1) f(-t)\,\overline{g(t-u)}\,K_{A_2}^j(t_2,w_2)dt\\
&\qquad= \int_{\mathbb R^2} K_{A_1}^i(-x_1,w_1) f(x)\,\overline{g(-x-u)}\,K_{A_2}^j(-x_2,w_2)(-dx)\\
&\qquad= \int_{\mathbb R^2} K_{A_1}^i(-x_1,w_1) f(x)\,\overline{g(-x-u)}\,K_{A_2}^j(-x_2,w_2)(-dx)\\
&\qquad=\dfrac{1}{2\pi \sqrt{b_1b_2}}\int_{\mathbb R^2} \,e^{\frac{i}{2b_1}\big[a_1(-x_1)^2-2(-x_1)(w_1-p_1)-2w_1(d_1p_1-b_1q_1)+d_1(w_1^2+p_1^2)-\frac{\pi b_1}{2}\big]}\\
&\qquad\qquad\times f(x)\,e^{\frac{j}{2b_2}\big[a_2(-x_2)^2-2(-x_2+k_2)(w_2-p_2)-2w_2(d_2p_2-b_2q_2)+d_2(w_2^2+p_2^2)-\frac{\pi b_2}{2}\big]}\,(-dx)\\
&\qquad=\dfrac{(-1)}{2\pi \sqrt{b_1b_2}}\int_{\mathbb R^2} \,e^{\frac{i}{2b_1}\big[a_1x_1^2-2(x_1)\big(-w_1-(-p_1)\big)-2(-w_1)\big(d_1(-p_1)-b_1(-q_1)\big)+d_1\big(w_1^2+(-p_1)^2\big)\big]}\\
&\qquad\qquad\qquad\qquad\times e^{-i\frac{\pi}{4}}\,f(x)\overline{g\big(-(x+u)\big)}e^{-j\frac{\pi}{4}}\\
&\qquad\quad\times \,e^{\frac{j}{2b_2}\big[a_2x_2^2-2(x_2)\big(-w_2-(-p_2)\big)-2(-w_2)\big(d_2(-p_2)-b_2(-q_2)\big)+d_2\big(w_2^2+(-p_2)^2\big)\big]}\,dx\\
&\qquad= -\int_{\mathbb R^2} K_{A^\prime_1}^i(x_1,-w_1) f(x)\,\overline{Pg\big(x-(-u)\big)}\,K_{A^\prime_2}^j(x_2,-w_2)\,dx\\
&\qquad= -{\mathcal{O}}^{A_1^\prime,A_2^\prime}_{Pg}\big[f(t)\big](-u,-w),\,\, \text{where}\quad A_s^\prime =\left[\begin{array}{cccc}a_s & b_s &| & -p_s\\c_s & d_s &| & -q_s \\\end{array}\right],\,s=1,2.
\end{align*}

\parindent=0mm \vspace{.1in}
In particular, it can easily be inferred that
\begin{align*}
{\mathcal{O}}^{A_1,A_2}_{Pg}\Big[P f(t)\Big](u,w)= -{\mathcal{O}}^{A_1^\prime,A_2^\prime}_{g}\big[f(t)\big](-u,-w).
\end{align*}

\parindent=0mm \vspace{.1in}
{\bf Property 5.} {\it (Anti-Linearity) For given quaternion window functions $g_n \in L^2(\mathbb R^2,\mathbb H), n\in \mathbb N,$ and $f \in L^2(\mathbb R^2,\mathbb H)$, we have:}
\begin{align*}
\mathcal O_{\sum \alpha_n g_n}^{A_1,A_2}[f](u,w)= \sum_{n\in \mathbb N}\, \mathcal O_{g_n}^{A_1,A_2}[f](u,w) \,\overline{\alpha}_n,\,\alpha_n\in \mathbb H.
\end{align*}

\parindent=0mm \vspace{.1in}
{\bf Proof.} This property follows similarly as the Property 1.

\parindent=0mm \vspace{.1in}
{\bf Property 6.} {\it (Phase-Shift) For any given quaternion function $f \in L^2(\mathbb R^2,\mathbb H),$ and window function $g \in L^2(\mathbb R^2,\mathbb H),$ we have:}
\begin{align*}
{\mathcal{O}}^{A_1,A_2}_{g(t-k)}[f](u,w)&={\mathcal{O}}^{A_1,A_2}_g  \big[f(t)\big](u+k,w)
\end{align*}
{\it where $u=(u_1,u_2),\,k=(k_1,k_2),\,w=(w_1,w_2).$}

\parindent=0mm \vspace{.1in}
{\bf Proof.} From Definition 3.1, we have
\begin{align*}
{\mathcal{O}}^{A_1,A_2}_{g(t-k)}[f](u,w)&= \int_{\mathbb R^2} K_{A_1}^i(t_1,w_1) f(t)\,\overline{g(t-k-u)}\,K_{A_2}^j(t_2,w_2)dt\\
&= \int_{\mathbb R^2} K_{A_1}^i(t_1,w_1) f(t)\,\overline{g\big(t-(u+k)\big)}\,K_{A_2}^j(t_2,w_2)dt\\
&={\mathcal{O}}^{A_1,A_2}_g  \big[f(t)\big](u+k,w)
\end{align*}

In the remaining part of this section, we will establish an inner product relation between two quaternion valued signals and their respective quaternion windowed OLCTs. As a consequence of this relation, we can deduce the Plancheral identity for the quaternion windowed OLCT as defined by (3.1).

\parindent=0mm \vspace{.1in}
{\bf Theorem 3.2( Inner Product Relation).} {\it Let ${\mathcal{O}}^{A_1,A_2}_{g_1} [f_1]$  and ${\mathcal{O}}^{A_1,A_2}_{g_2} [f_2]$ be the quaternion windowed OLCTs of 2D quaternion-valued signals $f_{1}$ and $f_{2}$, respectively. Then, we have}
\begin{align*}
\Big\langle {\mathcal{O}}^{A_1,A_2}_{g_1} [f_1],\,{\mathcal{O}}^{A_1,A_2}_{g_2} [f_2]\Big\rangle_{L^2(\mathbb R^2\times\mathbb R^2,\mathbb H)}=\big[\big\langle f_1,f_2 \big\rangle\overline{\big\langle g_1,g_2 \big\rangle}\big]_{\mathbb H}.\tag{3.2}
\end{align*}

\parindent=0mm \vspace{.0in}
{\bf Proof.} From definition 3.1, we have
\begin{align*}
&\Big\langle {\mathcal{O}}^{A_1,A_2}_{g_1} [f_1],\,{\mathcal{O}}^{A_1,A_2}_{g_2} [f_2]\Big\rangle_{L^2(\mathbb R^2\times\mathbb R^2,\mathbb H)}\\
&\qquad\quad=\int_{\mathbb R^2}\int_{\mathbb R^2} \left[{\mathcal{O}}^{A_1,A_2}_{g_1} [f_1] (u,w)\, \overline{{\mathcal{O}}^{A_1,A_2}_{g_2} [f_2](u,w)}\,\right]_{\mathbb H}\,du\,dw\\
&\qquad\quad=\int_{\mathbb R^4} \left[{\mathcal{O}}^{A_1,A_2}_{g_1} [f_1] (u,w)\,\int_{\mathbb R^2} \overline{ K_{A_1}^i(x_1,w_1) \,f_2(x)\overline{g_2(x-u)}\,K_{A_2}^j(x_2,w_2)dx}\,\right]_{\mathbb H}\,du\,dw\\
&\qquad\quad=\int_{\mathbb R^6} \left[{\mathcal{O}}^{A_1,A_2}_{g_1} [f_1] (u,w)\, \overline{K_{A_2}^j(x_2,w_2)} \,g_2(x-u)\,\overline{f_2(x)}\,\overline{ K_{A_1}^i(x_1,w_1)}\right]_{\mathbb H}dx\,du\,dw\\
&\qquad\quad=\int_{\mathbb R^8} \left[ K_{A_1}^i(y_1,w_1) \,f_1(y)\overline{g_1(y-u)}\,K_{A_2}^j(y_2,w_2) \overline{K_{A_2}^j(x_2,w_2)} \,g_2(x-u)\,\overline{f_2(x)}\,\overline{ K_{A_1}^i(x_1,w_1)}\,\right]_{\mathbb H}\\
&\qquad\qquad\qquad\qquad\qquad\qquad\qquad\qquad\qquad\qquad\qquad\qquad\qquad\qquad\qquad\qquad\times\,dy\,dx\,du\,dw\\
&\qquad\quad=\int_{\mathbb R^6} \left[ \left(\int_{\mathbb R^2}\overline{ K_{A_1}^i(x_1,w_1)} K_{A_1}^i(y_1,w_1)dw_1\right) \,f_1(y)\overline{g_1(y-u)}\right.\\
&\qquad\qquad\qquad\qquad\qquad\qquad\times\left.\left(\int_{\mathbb R^2}K_{A_2}^j(y_2,w_2) \overline{K_{A_2}^j(x_2,w_2)}dw_2 \right)\,g_2(x-u)\,\overline{f_2(x)}\,\right]_{\mathbb H}dxdydu\\
&\qquad\quad=\int_{\mathbb R^6} \left[ \delta(x_1-y_1) \,f_1(y)\overline{g_1(y-u)}\, \delta(x_2-y_2)\,g_2(x-u)\,\overline{f_2(x)}\,\right]_{\mathbb H}dxdydu\\
&\qquad\quad=\int_{\mathbb R^4} \left[ \,f_1(x)\overline{g_1(x-u)}\,g_2(x-u)\,\overline{f_2(x)}\,\right]_{\mathbb H}dxdu\\
&\qquad\quad=\int_{\mathbb R^4} \left[ \overline{f_2(x)}\,f_1(x)\overline{g_1(x-u)}\,g_2(x-u)\,\right]_{\mathbb H}dxdu\\
&\qquad\quad= \left[ \int_{\mathbb R^2}\overline{f_2(x)}\,f_1(x)dx\,\int_{\mathbb R^2}\overline{g_1(x-u)}\,g_2(x-u)\,du\,\right]_{\mathbb H}\\
&\qquad\quad= \left[\big\langle f_1,f_2 \big\rangle\overline{\big\langle g_1,g_2 \big\rangle}\right]_{\mathbb H}.
\end{align*}

This completes the proof of Theorem 3.2.\quad \fbox

\parindent=8mm \vspace{.1in}

The following corollary follows directly from Theorem 3.2.

\parindent=0mm \vspace{.1in}

{\bf Corollary 3.3(Energy Conservation).} {\it If $f,\,g\in L^2(\mathbb R^2,\mathbb H)$, then}
\begin{align*}
\int_{\mathbb R^2}\int_{\mathbb R^2}\Big|\left[{\mathcal{O}}^{A_1,A_2}_g f\right](u,w)\Big|^2 du\,dw=\big\|f\big\|^2_{L^2(\mathbb R^2,\mathbb H)}\big\|g\big\|^2_{L^2(\mathbb R^2,\mathbb H)},\tag{3.3}
\end{align*}

\parindent=0mm \vspace{.1in}
{\it Remarks:} For any normalized window function $g\in L^2(\mathbb R^2,\mathbb H)$,  the quaternion windowed OLCT becomes an isometry from the space of signals to the space of transforms or more precisely from $L^2(\mathbb R^2,\mathbb H)$ to $L^2(\mathbb R^4,\mathbb H)$.

\parindent=8mm \vspace{.1in}
Our next theorem guarantees the reconstruction of the input signal from the corresponding quaternion windowed OLCT.

\parindent=0mm \vspace{.1in}
{\bf Theorem 3.4. (Inversion Formula).} {\it Given any two quaternion-valued window functions $g_{1}$ and $g_{2}$ such that $\big\langle g_1,g_2 \big\rangle\ne 0$. Then, any $f\in L^2(\mathbb R^2,\mathbb H)$ can be reconstructed by the formula:}
\begin{align*}
f(t)=\dfrac{1}{\big\langle g_2,g_1 \big\rangle}\int_{\mathbb R^2}\int_{\mathbb R^2} \Big[K_{A_1}^{-i}(t_1,w_1){\mathcal{O}}^{A_1,A_2}_{g_1} [f](u,w)\,\overline{{\mathcal K}_{A_2}(t_2,w_2)}\cdot g_{2}(t-u)\Big]_{\mathbb H}\,dw\,du.\tag{3.4}
\end{align*}

\parindent=0mm \vspace{.0in}
{\it Proof.} As ${\mathcal{O}}^{A_1,A_2}_{g_1} [f]\in L^2(\mathbb R^4,\,\mathbb H)$, so we can easily deduce that the integral
\begin{align*}
\tilde f(t)=\dfrac{1}{\big\langle g_2,g_1 \big\rangle}\int_{\mathbb R^2}\int_{\mathbb R^2}\left[K_{A_1}^{-i}(t_1,w_1){\mathcal{O}}^{A_1,A_2}_{g_1} [f](u,w)\,K_{A_2}^{-j}(t_2,w_2)\, g_{2}(t-u)\right]_{\mathbb H}\,dw\,du
\end{align*}
is well defined in $L^2(\mathbb R^2,\,\mathbb H)$. Therefore, by virtue of the orthogonality relation (3.2), we observe that

\begin{align*}
\big\langle \tilde  f,\,h \big\rangle &=\dfrac{1}{\big\langle g_2,g_1 \big\rangle}\int_{\mathbb R^2}\int_{\mathbb R^2}\int_{\mathbb R^2}\left[K_{A_1}^{-i}(t_1,w_1){\mathcal{O}}^{A_1,A_2}_{g_1} [f](u,w)\,K_{A_2}^{-j}(t_2,w_2)\, g_{2}(t-u)\,\overline{h(x)}\,\right]_{\mathbb H}\,dw\,du\, dx\\
&=\dfrac{1}{\big\langle g_2,g_1 \big\rangle}\int_{\mathbb R^2}\int_{\mathbb R^2}\int_{\mathbb R^2}\left[{\mathcal{O}}^{A_1,A_2}_{g_1} [f](u,w)\,K_{A_2}^{-j}(t_2,w_2)\, g_{2}(t-u)\,\overline{h(x)}\,K_{A_1}^{-i}(t_1,w_1)\,\right]_{\mathbb H}\,dw\,du\, dx\\
&=\dfrac{1}{\big\langle g_2,g_1 \big\rangle}\int_{\mathbb R^2}\int_{\mathbb R^2}\left[{\mathcal{O}}^{A_1,A_2}_{g_1} [f](u,w)\,\int_{\mathbb R^2}\,K_{A_2}^{-j}(t_2,w_2)\, g_{2}(t-u)\,\overline{h(x)}\,K_{A_1}^{-i}(t_1,w_1)\, dx\,\right]_{\mathbb H}\,dw\,du\\
&=\dfrac{1}{\big\langle g_2,g_1 \big\rangle}\int_{\mathbb R^2}\int_{\mathbb R^2}\left[{\mathcal{O}}^{A_1,A_2}_{g_1} [f](u,w)\,\int_{\mathbb R^2}\overline{K_{A_1}^{i}(t_1,w_1)\,h(x)\,\overline{g_{2}(t-u)}\,K_{A_2}^{j}(t_2,w_2)}\, dx\,\right]_{\mathbb H}\,dw\,du\\
&=\dfrac{1}{\big\langle g_2,g_1 \big\rangle}\int_{\mathbb R^2}\int_{\mathbb R^2}\left[{\mathcal{O}}^{A_1,A_2}_{g_1} [f](u,w)\,\overline{{\mathcal{O}}^{A_1,A_2}_{g_2} [h](u,w)}\,\right]_{\mathbb H}\,dw\,du\\
&=\dfrac{1}{\big\langle g_2,g_1 \big\rangle}\left[\big\langle f_1,f_2 \big\rangle\overline{\big\langle g_1,g_2 \big\rangle}\right]_{\mathbb H}\\
&=\langle f_1,f_2 \big\rangle.
\end{align*}
Thus, $\tilde f=f$. This completes the proof of Theorem 3.4 \qquad \fbox

\parindent=0mm \vspace{.1in}

{\bf Corollary 3.5.} If we take $g_{1}=g_{2}=g$, then equation (3.4) reduces to
\begin{align*}
f(t)=\dfrac{1}{\| g\|^2} \int_{\mathbb R^2}\int_{\mathbb R^2} \Big[K_{A_1}^{-i}(t_1,w_1){\mathcal{O}}^{A_1,A_2}_{g} [f](u,w)\,{\mathcal K}^{-j}_{A_2}(t_2,w_2)\cdot g(t-u)\Big]_{\mathbb H}\,dw\,du.\tag{3.5}
\end{align*}

\section{Uncertainty Principles for Quaternion WOLCT}

\parindent=0mm \vspace{.1in}
The uncertainty principles in harmonic analysis are of central importance as they provide a lower bound for optimal simultaneous resolution in the time and frequency domains (see \cite{Fol}). The most famous of them is Heisenberg-Pauli-Weyl inequality, this has been extended to different time-frequency transforms and several other versions of the uncertainty principle have been investigated from time to time. For instance, Beckner \cite{Bec} obtained a logarithmic version of the uncertainty principle by using a sharp form of Pitt's inequality and showed that this version yields the classical Heisenberg's inequality by virtue of Jensen's inequality. Local version of uncertainty states that if a function is concentrated, then not only is its transformation spread out, but that cannot be localized in a subset of finite measure (see \cite{Lian}).  In this Section, we shall establish an analogue of the well-known Heisenberg's uncertainty inequality and the corresponding logarithmic version and the local uncertainty principle for the quaternion windowed offset linear canonical transform as defined by (3.1). First, we prove the following lemma

\parindent=0mm \vspace{.1in}

{\bf Lemma 4.1.} {\it Let $g\in L^2(\mathbb R^2,\,\mathbb H)\backslash \{0\}$ be a quaternion window function, then for every $f \in L^2(\mathbb R^2,\,\mathbb H),$ we have }

\begin{align*}
\|g\|^2_{L^2(\mathbb R^2,\,\mathbb H)} \, \int_{\mathbb R^2} t_k^2\,\big|f(t)\big|_{\mathbb H}^2\, dt = \int_{\mathbb R^2} t_k^2 \,\Big|\mathcal{O}^{-1}_{A_1,A_2} {\Big[\mathcal O_g^{A_1,A_2} [f](u,\xi)\Big]}(t)\Big|_{\mathbb H}^2  dt\,du\,\tag{4.1}
\end{align*}

\parindent=0mm \vspace{.0in}
{\it Proof.} Using Corollary 3.5, we observe that
\begin{align*}
&\left\|g\right\|^2_{L^2(\mathbb R^2,\,\mathbb H)} \,\int_{\mathbb R^2} t_k^2\big|f(t)\big|_{\mathbb H}^2 dt
\\&\qquad= \int_{\mathbb R^2} t_k^2\big|f(t)\big|_{\mathbb H}^2 dt \int_{\mathbb R^2}\big|g(t-b)\big|_{\mathbb H}^2 du
\\&\qquad=\int_{\mathbb R^2}\int_{\mathbb R^2} t_k^2\big|f(t)\big|_{\mathbb H}^2 \big|g(t-b)\big|_{\mathbb H}^2 dt\,du
\\&\qquad=\int_{\mathbb R^2}\int_{\mathbb R^2} t_k^2\big|f(t)\,\overline{g(t-b)}\big|_{\mathbb H}^2 dt\,du
\\&\qquad=\int_{\mathbb R^4} t_k^2\left|\dfrac{1}{\| g\|^2} \int_{\mathbb R^4} \Big[K_{A_1}^{-i}(t_1,\xi_1){\mathcal{O}}^{A_1,A_2}_{g} [f](u,\xi)\,{\mathcal K}^{-j}_{A_2}(t_2,\xi_2)\cdot g(t-u)\overline{g(t-u)}\Big]_{\mathbb H}\,d\xi\,du\right|_{\mathbb H}^2 dtdb
\\&\qquad=\int_{\mathbb R^4} t_k^2\left|\dfrac{1}{\| g\|^2} \int_{\mathbb R^2} \Big[K_{A_1}^{-i}(t_1,\xi_1){\mathcal{O}}^{A_1,A_2}_{g} [f](u,\xi)\,{\mathcal K}^{-j}_{A_2}(t_2,\xi_2)\int_{\mathbb R^2}|g(t-u)|_{\mathbb H}^2\,du\,\Big]_{\mathbb H}\,d\xi\,\right|_{\mathbb H}^2 dt\,db
\\&\qquad=\int_{\mathbb R^4} t_k^2\left| \int_{\mathbb R^2} \Big[K_{A_1}^{-i}(t_1,\xi_1){\mathcal{O}}^{A_1,A_2}_{g} [f](u,\xi)\,{\mathcal K}^{-j}_{A_2}(t_2,\xi_2)\Big]_{\mathbb H}\,d\xi\,\right|_{\mathbb H}^2 dt\,db
\\&\qquad= \int_{\mathbb R^2} t_k^2 \,\Big|\mathcal{O}^{-1}_{A_1,A_2} {\Big[\mathcal O_g^{A_1,A_2} [f](u,\xi)\Big]}(t)\Big|_{\mathbb H}^2  dt\,du
\end{align*}
This completes the proof of Lemma 4.1.\quad \fbox

\parindent=8mm \vspace{.2in}
We are now ready to establish the Heisenberg-type inequalities for the proposed quaternion windowed offset linear canonical transform as defined by (3.1).

\parindent=0mm \vspace{.1in}
{\bf Theorem 4.2.} {\it Let $g\in L^2(\mathbb R^2,\,\mathbb H)\backslash \{0\}$ be a quaternion window function $\mathcal{O}^{A_1,A_2}_g [f]$ be the quaternion windowed offset linear canonical transform of any signal $f \in L^2(\mathbb R^2,\,\mathbb H)$, then we have:}
\begin{align*}
\left\{\int_{\mathbb R}t_k^2 \big|f(t)\big|_{\mathbb H}^2 dt\right\}^{1/2} \left\{\int_{\mathbb R^2} \int_{\mathbb R^2}\xi_k^2 \Big|\mathcal O_g^{A_1,A_2} [f](u,\xi)\Big|_{\mathbb H}^2 dud\xi \right\}^{1/2} \ge \dfrac{b_k}{2}\big\|f\big\|^2_{L^2(\mathbb R^2,\,\mathbb H)}\big\|g\big\|^2_{L^2(\mathbb R^2,\,\mathbb H)} \tag{4.2}
\end{align*}

\parindent=0mm \vspace{.1in}
{\it Proof.} By virtue of the Heisenberg's  inequality for the quaternion offset linear canonical transform \cite{QOLCT}, we can write
\begin{align*}
\int_{\mathbb R^2} t^2\big|f(t)\big|_{\mathbb H}^2dt \int_{\mathbb R^2} \xi_k^2 \Big|\mathcal{O}^{A_1,A_2} [f](\xi)\Big|_{\mathbb H}^2 d\xi\ge \left\{\dfrac{b_k}{2} \int_{\mathbb R}\big|f(t)\big|_{\mathbb H}^2dt \right\}^2
\end{align*}

By implementing the inversion formula and the orthogonality relation for quaternion offset linear canonical transform, we can rewrite the above inequality as

\begin{align*}
\int_{\mathbb R^2} t^2\left|\mathcal{O}^{-1}_{A_1,A_2} \Big[\mathcal{O}^{A_1,A_2} [f](\xi)\Big](t)\right|_{\mathbb H}^2dt \int_{\mathbb R^2} \xi_k^2 \Big|\left[\mathcal{O}^{A_1,A_2} f\right](\xi)\Big|_{\mathbb H}^2 d\xi \ge \left\{\dfrac{b_k}{2} \int_{\mathbb R^2}\Big|\mathcal{O}^{A_1,A_2} [f](\xi)\Big|_{\mathbb H}^2d\xi \right\}^2
\end{align*}

Since  $\mathcal{O}^M_g \left[f\right]\in L^2(\mathbb R,\,\mathbb H)$, therefore we can replace $\mathcal{O}^M [f]$ by $\mathcal{O}^M_g [f]$ to obtain

\begin{align*}
&\int_{\mathbb R^2} t^2\left|\mathcal{O}^{-1}_{A_1,A_2} \Big[\mathcal{O}^{A_1,A_2}_{g} [f](u,\xi)\Big](t)\right|_{\mathbb H}^2dt \int_{\mathbb R^2} \xi_k^2 \Big|\left[\mathcal{O}^{A_1,A_2}_{g} f\right](\xi)\Big|_{\mathbb H}^2 d\xi\\
&\qquad\qquad\qquad\qquad\qquad\qquad\qquad\qquad\qquad\qquad\qquad\qquad\ge \left\{\dfrac{b_k}{2} \int_{\mathbb R^2}\Big|\mathcal{O}^{A_1,A_2}_{g} [f](u,\xi)\Big|_{\mathbb H}^2d\xi \right\}^2
\end{align*}

Taking square root on both sides of above equation and integrating with respect to measure $du$, we obtain

\begin{align*}
&\int_{\mathbb R^2} \Bigg(\left\{\int_{\mathbb R^2}t^2\left|\mathcal{O}^{-1}_{A_1,A_2} \Big[\mathcal{O}^{A_1,A_2}_{g} [f](u,\xi)\Big](t)\right|_{\mathbb H}^2dt\right\}^{1/2} \left\{\int_{\mathbb R^2} \xi_k^2 \Big|\left[\mathcal{O}^{A_1,A_2}_{g} f\right](\xi)\Big|_{\mathbb H}^2 d\xi\,du\right\}^{1/2}\Bigg)\,du\\
&\qquad\qquad\qquad\qquad\qquad\qquad\qquad\qquad\qquad\qquad\qquad\ge \dfrac{b_k}{2} \int_{\mathbb R^2}\int_{\mathbb R^2}\Big|\mathcal{O}^{A_1,A_2}_{g} [f](u,\xi)\Big|_{\mathbb H}^2d\xi\,du
\end{align*}

Furthermore, an implication of the well known Cauchy-Schwarz inequality yields

\begin{align*}
&\left\{\int_{\mathbb R^2} \int_{\mathbb R^2}t_k^2\left|\mathcal{O}^{-1}_{A_1,A_2} \Big[\mathcal{O}^{A_1,A_2}_{g} [f](u,\xi)\Big](t)\right|_{\mathbb H}^2dt\,du\right\}^{1/2} \left\{\int_{\mathbb R^2}\int_{\mathbb R^2} \xi_k^2 \Big|\left[\mathcal{O}^{A_1,A_2}_{g} f\right](\xi)\Big|_{\mathbb H}^2 d\xi\,du\right\}^{1/2}\\
&\qquad\qquad\qquad\qquad\qquad\qquad\qquad\qquad\qquad\qquad\qquad\qquad\ge \dfrac{b_k}{2} \int_{\mathbb R^2}\int_{\mathbb R^2}\Big|\mathcal{O}^{A_1,A_2}_{g} [f](u,\xi)\Big|_{\mathbb H}^2d\xi\,du
\end{align*}

Finally, on implementing the Lemma 4.1 on L.H.S, and Corollary 3.3 on R.H.S of the above inequality we have

\begin{align*}
&\left\{\left\|g\right\|^2_{L^2(\mathbb R^2,\,\mathbb H)} \,\int_{\mathbb R^2} t_k^2\big|f(t)\big|_{\mathbb H}^2 dt \right\}^{1/2} \left\{\int_{\mathbb R^2}\int_{\mathbb R^2} \xi_k^2 \Big|\left[\mathcal{O}^{A_1,A_2}_{g} f\right](\xi)\Big|_{\mathbb H}^2 d\xi\,du\right\}^{1/2}\qquad\\
&\qquad\qquad\qquad\qquad\qquad\qquad\qquad\qquad\qquad\qquad\qquad\qquad\ge \dfrac{b_k}{2}\big\|f\big\|^2_{L^2(\mathbb R^2,\mathbb H)}\big\|g\big\|^2_{L^2(\mathbb R^2,\mathbb H)}
\end{align*}

The desired result is obtained by dividing $\big\|g\big\|^2_{L^2(\mathbb R^2,\mathbb H)}$ on both sides of above inequality. \qquad\fbox

\parindent=8mm \vspace{.1in}

We now establish the logarithmic uncertainty principle for the quaternion windowed offset linear canonical transform $\mathcal{O}^M_g [f]$ as defined by (3.1). First, we remind the following definition of space of rapidly decreasing smooth quaternion functions(see \cite{Asym}).

\parindent=0mm \vspace{.1in}
{\bf Definition 4.3.} For a multi-index $\alpha =(\alpha_1,\alpha_2)\in \mathbb R^{+}\times \mathbb R^{+},$  the Schwartz space in $ L^2(\mathbb R^2,\mathbb H)$ is defined as

\begin{align*}
\mathcal S (\mathbb R^2,\mathbb H)= \left\{f\in\mathbb C^{\infty} (\mathbb R^2,\mathbb H); \sup_{t\in\mathbb R^2}\big(1+|t|^k\big)\Big|\dfrac{\partial^{\alpha_1+\alpha_2}[f(t)]}{\partial^{\alpha_1}_{t_1}\partial^{\alpha_2}_{t_2}}\Big|\,< \infty \right\},
\end{align*}
where $\mathbb C^{\infty}(\mathbb R^2,\mathbb H)$ is the set of smooth functions from $\mathbb R^2$ to $\mathbb H$.

\parindent=0mm \vspace{.1in}

{\bf Theorem 4.4.}  {\it Given a quaternion window function $g\in \mathcal S (\mathbb R^2,\mathbb H)$ and a signal $f\in \mathcal S(\mathbb R^2,\mathbb H)$, the quaternion windowed offset linear canonical transform $\mathcal{O}^{A_1,A_2}_g [f](u,\xi)$ satisfies the following logarithmic estimate of the uncertainty inequality:}
\begin{align*}
&\int_{\mathbb R^{2}}\int_{\mathbb R^{2}}\ln\Big|\frac{\xi}{b}\Big|\, \Big|\mathcal{O}^{A_1,A_2}_g [f](u,\xi)\Big|_{\mathbb H}^{2} dud\xi +\|g\|^2_{L^2(\mathbb R^2,\mathbb H)} \int_{\mathbb R^2}\ln|t|\,|f(t)|_{\mathbb H}^{2}dt\qquad\qquad\qquad\qquad\qquad\\
&\qquad\qquad\qquad\qquad\qquad\qquad\qquad\qquad\qquad\qquad\qquad\qquad\geq\,\mathcal D\,\|g\|^2_{L^2(\mathbb R^2,\mathbb H)}\,\|f\|^2_{L^2(\mathbb R^2,\mathbb H)} \tag{4.3}
\end{align*}
{\it where $\mathcal D=\left(\dfrac{\Gamma^{\prime}(1/2)}{\Gamma(1/2)}-\ln 2\right)$ and $\Gamma$ is a Gamma function.}

\parindent=0mm \vspace{.1in}
{\it Proof.} For the quaternion-valued function $f\in L^2(\mathbb R^2,\mathbb H)$, the time and frequency spreads satisfy the inequality \cite{QOLCT}
\begin{align*}
\int_{\mathbb R^{2}}\ln|t|\,\big|f(t)\big|^{2}dt+\int_{\mathbb R^{2}}\ln\Big|\frac{\xi}{b}\Big|\,{\Big|\mathcal{O}^{A_1,A_2} [f](u,\xi)\Big|_{\mathbb H}^{2}}\,d\xi\geq\,\mathcal D\int_{\mathbb R^{2}}\Big|f(t)\Big|_{\mathbb H}^{2}\,dt.
\end{align*}

Invoking the inversion formula of QOLCT on the L.H.S and Parseval's formula for QOLCT on R.H.S, we obtain
\begin{align*}
\int_{\mathbb R^{2}}\ln|t|\left|\mathcal{O}^{-1}_{A_1,A_2} \Big[\mathcal{O}^{A_1,A_2} [f](\xi)\Big](t)\right|_{\mathbb H}^{2} dt+\int_{\mathbb R^{2}}\ln\Big|\frac{\xi}{b}\Big|{\Big|\mathcal{O}^{A_1,A_2} [f](\xi)\Big|_{\mathbb H}^{2}}d\xi\geq\mathcal D\int_{\mathbb R^{2}}\Big|\mathcal O^{A_1,A_2}[f](\xi)\Big|_{\mathbb H}^{2}d\xi
\end{align*}

Replacing $\mathcal O^{A_1,A_2}[f](\xi)$ by  $\mathcal{O}^{A_1,A_2}_g [f](u,\xi)$ in the above inequality, we obtain
\begin{align*}
&\int_{\mathbb R^{2}}\ln|t|\,\left|\mathcal{O}^{-1}_{A_1,A_2} \Big[\mathcal{O}^{A_1,A_2}_g [f](u,\xi)\Big](t)\right|_{\mathbb H}^{2} dt+\int_{\mathbb R^{2}}\ln\Big|\frac{\xi}{b}\Big|\,{\Big|\mathcal{O}^{A_1,A_2}_g [f](u,\xi)\Big|_{\mathbb H}^{2}}\,d\xi\qquad\qquad\qquad\qquad\qquad\\
&\qquad\qquad\qquad\qquad\qquad\qquad\qquad\qquad\qquad\qquad\qquad\qquad\geq\,\mathcal D\int_{\mathbb R^{2}}\Big|\mathcal{O}^{A_1,A_2}_g [f](u,\xi)\Big|_{\mathbb H}^{2}\,d\xi.
\end{align*}
Integrating above equation with respect to measure $du$ and then applying the  Fubini theorem, we obtain
\begin{align*}
&\int_{\mathbb R^{2}}\int_{\mathbb R^{2}}\ln|t|\,\left|\mathcal{O}^{-1}_{A_1,A_2} \Big[\mathcal{O}^{A_1,A_2}_g [f](u,\xi)\Big](t)\right|_{\mathbb H}^{2} dtdu+\int_{\mathbb R^{2}}\int_{\mathbb R^{2}}\ln\Big|\frac{\xi}{b}\Big|\,{\Big|\mathcal{O}^{A_1,A_2}_g [f](u,\xi)\Big|_{\mathbb H}^{2}}\,d\xi\,du\qquad\\
&\qquad\qquad\qquad\qquad\qquad\qquad\qquad\qquad\qquad\geq\,\mathcal D\int_{\mathbb R^{2}}\int_{\mathbb R^{2}}\Big|\mathcal{O}^{A_1,A_2}_g [f](u,\xi)\Big|_{\mathbb H}^{2}\,d\xi\,du.
\end{align*}

Applying Lemma 4.1 on L.H.S and Corollary 3.6 on R.H.S, we obtain the desired result
\begin{align*}
&\int_{\mathbb R^{2}}\int_{\mathbb R^{2}}\ln\Big|\frac{\xi}{b}\Big|\, \Big|\mathcal{O}^{A_1,A_2}_g [f](u,\xi)\Big|_{\mathbb H}^{2} dud\xi +\|g\|^2_{L^2(\mathbb R^2,\mathbb H)} \int_{\mathbb R^2}\ln|t|\,|f(t)|_{\mathbb H}^{2}dt\qquad\qquad\qquad\qquad\qquad\\
&\qquad\qquad\qquad\qquad\qquad\qquad\qquad\qquad\qquad\qquad\qquad\qquad\geq\,\mathcal D\,\|g\|^2_{L^2(\mathbb R^2,\mathbb H)}\,\|f\|^2_{L^2(\mathbb R^2,\mathbb H)}
\end{align*}

This completes the proof of Theorem 4.4. \quad\fbox

\parindent=8mm \vspace{.1in}

In the following, we establish a local type uncertainty principle for quaternion windowed offset linear canonical transform.

\parindent=0mm \vspace{.1in}

{\bf Theorem 4.5.}  {\it Given a quaternion window function $g\in L(\mathbb R^2,\mathbb H)$  and a signal
$f\in L(\mathbb R^2,\mathbb H)$, with $\|g\|^2_{L^2(\mathbb R^2,\mathbb H)}=1=\|f\|^2_{L^2(\mathbb R^2,\mathbb H)},$
such that for a measurable set $E\subset \mathbb R^2\times\mathbb R^2
,\, \epsilon\geq 0,$ and}
\begin{align*}
\int \int_{E} \Big|\mathcal{O}^{A_1,A_2}_g [f](u,\xi)\Big|_{\mathbb H}^{2} dud\xi\,\geq\,1-\epsilon.\tag{4.5}
\end{align*}
We have $2\pi (1-\epsilon)\big| b_1 b_2 \big|^{1/2}\leq \,\mu(E),$ where $\mu(E)$ is Lebesgue measure of $E$.

\parindent=0mm \vspace{.1in}
{\it Proof.} From Definition 3.1, we have
\begin{align*}
\Big|\mathcal O^{A_1,A_2}_{g}[f](u,\xi) \Big|_{\mathbb H}&=\Big| \int_{\mathbb R^2} K_{A_1}^i(t_1,\xi_1) \,f(t)\overline{g(t-u)}\,K_{A_2}^j(t_2,\xi_2)dt\Big|_{\mathbb H}
\\&\leq\,\dfrac{1}{2\pi \big| b_1 b_2 \big|^{1/2}}\int_{\mathbb R^2}|f(t)|\,\big|\overline{g(t-u)}\big|\,dx.
\end{align*}

Taking Sup-norm on L.H.S and implementing well known Holders inequality on R.H.S, we obtain

\begin{align*}
\Big\|\mathcal O^{A_1,A_2}_{g}[f](u,\xi) \Big\|_{L^{\infty}(\mathbb R^2,\,\mathbb H)}
\leq\,\dfrac{1}{2\pi \big| b_1 b_2 \big|^{1/2}}\,\big\|f\big\|_{L^{2}(\mathbb R^2,\,\mathbb H)}\big\|g\big\|_{L^{2}(\mathbb R^2,\,\mathbb H)}\,\tag{4.6}
\end{align*}

Plugging inequality (4.6) in (4.5), we get
\begin{align*}
1-\epsilon\,&\leq\,\int \int_{E} \Big|\mathcal{O}^{A_1,A_2}_g [f](u,\xi)\Big|_{\mathbb H}^{2} dud\xi\\
&\leq\, \mu(E)\cdot \Big\|\mathcal O^{A_1,A_2}_{g}[f](u,\xi) \Big\|_{L^{\infty}(\mathbb R^2,\,\mathbb H)}\\
&\leq\, \mu(E)\dfrac{1}{2\pi \big| b_1 b_2 \big|^{1/2}}\,\big\|f\big\|_{L^{2}(\mathbb R^2,\,\mathbb H)}\big\|g\big\|_{L^{2}(\mathbb R^2,\,\mathbb H)}\\
&\leq\, \mu(E)\dfrac{1}{2\pi \big| b_1 b_2 \big|^{1/2}}.
\end{align*}

This completes the proof of Theorem 4.5. \quad\fbox

\parindent=0mm \vspace{.1in}

{\bf Theorem 4.6(Local uncertainty inequality).}  {\it Let $E$ be a measurable subset of $\mathbb R^2\times\mathbb R^2$, such that $0<\mu(E)<1,$ then for every
$f\in L(\mathbb R^2,\mathbb H)$ and a quaternion window function $g\in L (\mathbb R^2,\mathbb H)$ , we have}

\begin{align*}
\big\|f\big\|_{L^{2}(\mathbb R,\,\mathbb H)}\big\|g\big\|_{L^{2}(\mathbb R,\,\mathbb H)}\leq \dfrac{1}{\sqrt{1-\mu(E)}}
\Big\|\mathcal{O}^{A_1,A_2}_g [f](u,\xi)\Big\|_{L^2(E^c,\,\mathbb H)}\,\tag{4.7}
\end{align*}

{\it Moreover, for every $\alpha>0,$ there exist $C(\alpha)>0,$ such that}

\begin{align*}
\big\|f\big\|_{L^{2}(\mathbb R,\,\mathbb H)}\big\|g\big\|_{L^{2}(\mathbb R,\,\mathbb H)}\leq C(\alpha)\,\Big[\int_{\mathbb R^2}\int_{\mathbb R^2}\big|(u,\xi)\big|^{2\alpha}\,
\Big|\mathcal{O}^{A_1,A_2}_g [f](u,\xi)\Big|^2_{\mathbb H}\,dud\xi\Big]^{1/2}.\,\tag{4.8}
\end{align*}

\parindent=0mm \vspace{.1in}
{\it Proof.} By invoking the Corollary 3.3 and Theorem 4.5, we have

\begin{align*}
\Big\|\mathcal{O}^{A_1,A_2}_g [f](u,\xi)\Big\|^2_{L^2(\mathbb R^2\times\mathbb R^2,\,\mathbb H)}
&=\int_{\mathbb R^2}\int_{\mathbb R^2} \Big|\mathcal{O}^{A_1,A_2}_g [f](u,\xi)\Big|^2_{\mathbb H}dud\xi\\
&=\int \int_{E} \Big|\mathcal{O}^{A_1,A_2}_g [f](u,\xi)\Big|^2_{\mathbb H}dud\xi+\int \int_{E^c} \Big|\mathcal{O}^{A_1,A_2}_g [f](u,\xi)\Big|^2_{\mathbb H}dud\xi\\
&\leq \,\mu(E) \big\|f\big\|_{L^{2}(\mathbb R,\,\mathbb H)}\big\|g\big\|_{L^{2}(\mathbb R,\,\mathbb H)}+\Big\|\mathcal{O}^{A_1,A_2}_g [f](u,\xi)\Big\|^2_{L^2(E^c,\,\mathbb H)}
\end{align*}

Equivalently,
\begin{align*}
\big\|f\big\|^2_{L^{2}(\mathbb R,\,\mathbb H)}\big\|g\big\|^2_{L^{2}(\mathbb R,\,\mathbb H)}(1-\mu(E))\leq \Big\|\mathcal{O}^{A_1,A_2}_g [f](u,\xi)\Big\|^2_{L^2(E^c,\,\mathbb H)}
\end{align*}

Taking square root on both sides and then dividing both sides by $\sqrt{1-\mu(E)},$ we get

\begin{align*}
\big\|f\big\|_{L^{2}(\mathbb R,\,\mathbb H)}\big\|g\big\|_{L^{2}(\mathbb R,\,\mathbb H)}\leq \dfrac{1}{\sqrt{1-\mu(E)}}\Big\|\mathcal{O}^{A_1,A_2}_g [f](u,\xi)\Big\|_{L^2(E^c,\,\mathbb H)}
\end{align*}

which proves our first assertion.

\parindent=0mm \vspace{.0in}
Now, we fix $r_0\in(0,1]$ small enough such that $\mu\big(B_{r_0}\big)<1,$ where
$$B_{r_0}=\left\{(u,\xi)\in \mathbb R^2\times\mathbb R^2;\,\big|(u,\xi)\big|<\,r_0 \right\},$$
the ball of radius $r_0$ centered at origin, we have from inequality (4.7),

\begin{align*}
\big\|f\big\|_{L^{2}(\mathbb R,\,\mathbb H)}\big\|g\big\|_{L^{2}(\mathbb R,\,\mathbb H)}
&\leq\, \dfrac{1}{\sqrt{1-\mu(B_{r_0})}}\left[\int \int_{\big|(u,\xi)\big|>\,r_0 }\Big|\mathcal{O}^{A_1,A_2}_g [f](u,\xi)\Big|^2_{\mathbb H}dud\xi\right]^{1/2}\\
&\leq\, \dfrac{1}{r^{\alpha}_0 \sqrt{1-\mu(B_{r_0})}}\left[\int \int_{\big|(u,\xi)\big|>\,r_0 }r^{2\alpha}_0\Big|\mathcal{O}^{A_1,A_2}_g [f](u,\xi)\Big|^2_{\mathbb H}dud\xi\right]^{1/2}\\
&\leq\, \dfrac{1}{r^{\alpha}_0 \sqrt{1-\mu(B_{r_0})}}\left[\int \int_{\big|(u,\xi)\big|>\,r_0} \,\big|(u,\xi)\big|^{2\alpha}\Big|\mathcal{O}^{A_1,A_2}_g [f](u,\xi)\Big|^2_{\mathbb H}dud\xi\right]^{1/2}\\
&\leq\, C(\alpha)\left[\int \int_{\big|(u,\xi)\big|>\,r_0} \,\big|(u,\xi)\big|^{2\alpha}\Big|\mathcal{O}^{A_1,A_2}_g [f](u,\xi)\Big|^2_{\mathbb H}dud\xi\right]^{1/2}
\end{align*}
where $C(\alpha)=\dfrac{1}{r^{\alpha}_{0} \sqrt{1-\mu(B_{r_0})}}.$ This completes the proof of  Theorem 4.6. \quad\fbox

\section{Example of Quaternion Windowed Offset Linear canonical Transform}

\parindent=0mm \vspace{.0in}
For the demonstration of the quaternion windowed offset linear canonical transform (3.1), in this section, we shall present an illustrative example.

\parindent=0mm \vspace{.1in}
{\bf Example 5.1.} Consider a two-dimensional quaternion valued Gaussian  function $$f(t)=\beta \exp\left\{-
(\alpha_1 t_1^2+\alpha_2 t_2^2)\right\}, \alpha_1,\alpha_2\,\in \mathbb C, \,\beta\in L^2(\mathbb R^2,\,\mathbb H)$$ Then, the quaternion windowed offset linear canonical transform of $f(t)$ with respect to
the unimodular matrices $A_s=\left(a_s,b_s,c_s,d_s, \,|\, p_s, q_s\right),\,s=1,2,$ and the rectangular window function
\begin{align*}
g(t)=\left\{\begin{array}{cc} 1,\, &\text{if}\,\, |t_1|<a, |t_2|<a, \,a>0,\\ 0,\, &\text{elsewhere}\\ \end{array}\right.
\end{align*}
is given by
\begin{align*}
&{\mathcal{O}}_g^{A_1,A_2} [f](u,w)\\
&\quad=\dfrac{1}{2\pi \sqrt{b_1b_2}}\int_{\mathbb R^2} \,e^{\frac{i}{2b_1}\big[a_1t_1^2-2t_1(w_1-p_1)-2w_1(d_1p_1-b_1q_1)+d_1(w_1^2+p_1^2)-\frac{\pi b_1}{2}\big]} \,\beta \,e^{-(\alpha_1 t_1^2+\alpha_2t_2^2)}\\
&\qquad\qquad\times\overline{g(t-u)} \,e^{\frac{j}{2b_2}\big[a_2t_2^2-2t_2(w_2-p_2)-2w_2(d_2p_2-b_2q_2)+d_2(w_2^2+p_2^2)-\frac{\pi b_2}{2}\big]}dt\\
&\quad=\dfrac{1}{2\pi \sqrt{b_1b_2}}\int_{u_1-a}^{u_1+a}\int_{u_2-a}^{u_2+a} \,e^{\frac{i}{2b_1}\big[a_1t_1^2-2t_1(w_1-p_1)-2w_1(d_1p_1-b_1q_1)+d_1(w_1^2+p_1^2)-\frac{\pi b_1}{2}\big]-\alpha_1 t_1^2} \\
&\qquad\qquad\times\,\beta \,e^{\frac{j}{2b_2}\big[a_2t_2^2-2t_2(w_2-p_2)-2w_2(d_2p_2-b_2q_2)+d_2(w_2^2+p_2^2)-\frac{\pi b_2}{2}\big]-\alpha_2t_2^2}dt\\
&\quad=\dfrac{1}{2\pi \sqrt{b_1b_2}}\,e^{\frac{i}{2b_1}\big[-2w_1(d_1p_1-b_1q_1)+d_1(w_1^2+p_1^2)-\frac{\pi b_1}{2}\big]}\\
&\quad\times\int_{u_1-a}^{u_1+a} \,e^{\frac{i}{2b_1}\big[t_1^2(a_1-2ib_1\alpha_1)-2t_1(w_1-p_1)\big] }dt_1 \beta \int_{u_2-a}^{u_2+a}\,e^{\frac{j}{2b_2}\big[t_2^2(a_2-2jb_2\alpha_2)-2t_2(w_2-p_2)\big]}\,dt_2  \\
&\qquad\qquad\times \,e^{\frac{j}{2b_2}\big[-2w_2(d_2p_2-b_2q_2)+d_2(w_2^2+p_2^2)-\frac{\pi b_2}{2}\big]}\\
\end{align*}
For simplicity, choose $\alpha_1=\frac{-ia_1}{2b_1}$ and $\alpha_2=\frac{-j a_2}{2b_2},$ we have
\begin{align*}
&{\mathcal{O}}_g^{A_1,A_2} [f](u,w)\\
&\quad=\dfrac{1}{2\pi \sqrt{b_1b_2}}\,e^{\frac{i}{2b_1}\big[-2w_1(d_1p_1-b_1q_1)+d_1(w_1^2+p_1^2)-\frac{\pi b_1}{2}\big]}\,\int_{u_1-a}^{u_1+a} \,e^{\frac{i(p_1-w_1)}{b_1}t_1 }\,dt_1   \\
&\qquad\qquad\qquad\times \,\beta \int_{u_2-a}^{u_2+a}\,e^{\frac{j(p_2-w_2)}{b_2}t_2}\,dt_2\,e^{\frac{j}{2b_2}\big[-2w_2(d_2p_2-b_2q_2)+d_2(w_2^2+p_2^2)-\frac{\pi b_2}{2}\big]}\\
&\quad=\dfrac{1}{2\pi \sqrt{b_1b_2}}\,e^{\frac{i}{2b_1}\big[2w_1(b_1q_1-d_1p_1)+d_1(w_1^2+p_1^2)-\frac{\pi b_1}{2}\big]}\,\dfrac{b_1}{i(p_1-w_1)}\left[ e^{\frac{i(p_1-w_1)(u_1+a)}{b_1}}-e^{\frac{i(p_1-w_1)(u_1-a)}{b_1}}\right]  \\
&\qquad\qquad\times \,\beta \,\dfrac{b_2}{j(p_2-w_2)}\left[ e^{\frac{j(p_2-w_2)(u_2+a)}{b_1}}-e^{\frac{j(p_2-w_2)(u_2-a)}{b_2}}\right] \,e^{\frac{j}{2b_2}\big[2w_2(b_2q_2-d_2p_2)+d_2(w_2^2+p_2^2)-\frac{\pi b_2}{2}\big]}\\
&\quad=\dfrac{\sqrt{b_1b_2}}{2\pi(w_1-p_1)}\,e^{\frac{i}{2b_1}\big[2w_1(b_1q_1-d_1p_1)+d_1(w_1^2+p_1^2)+\frac{\pi b_1}{2}\big]}\left[ e^{\frac{i(p_1-w_1)(u_1+a)}{b_1}}-e^{\frac{i(p_1-w_1)(u_1-a)}{b_1}} \right]  \\
&\qquad\times \dfrac{\beta}{(w_2-p_2)} \left[ e^{\frac{j(p_2-w_2)(u_2+a)}{b_1}}-e^{\frac{j(p_2-w_2)(u_2-a)}{b_2}} \right] \,e^{\frac{j}{2b_2}\big[2w_2(b_2q_2-d_2p_2)+d_2(w_2^2+p_2^2)+\frac{\pi b_2}{2}\big]}.\\
\end{align*}

\end{document}